\magnification 1200
\def\R{{\rm I\kern-0.2em R\kern0.2em \kern-0.2em}}
\def\N{{\rm I\kern-0.2em N\kern0.2em \kern-0.2em}}
\def\P{{\rm I\kern-0.2em P\kern0.2em \kern-0.2em}}
\def\B{{\rm I\kern-0.2em B\kern0.2em \kern-0.2em}}
\def\C{{\bf \rm C}\kern-.4em {\vrule height1.4ex width.08em depth-.04ex}\;}

\def\D{{\Delta}}

\font\ninerm=cmr8
\ 
\vskip 6mm
\centerline {\bf HOLOMORPHIC EXTENDIBILITY AND THE ARGUMENT PRINCIPLE}
\vskip 4mm
\centerline{Josip Globevnik}
\vskip 4mm
\centerline{\it Dedicated to the memory of Herb Alexander}
\vskip 4mm\vskip 4mm
{\noindent \ninerm ABSTRACT\ \  The paper gives the following 
characterization of the disc 
algebra in terms of the argument principle: A continuous function 
$f$ on the unit circle $b\D$
extends holomorphically 
through the unit disc if and only if for each polynomial $P$ such that $f+P\not= 0 $ on 
$b\D$ the change of argument of $f+P$ around $b\D $ is nonnegative. } 
\vskip 6mm
\bf 1.\ Introduction and the main result \rm 
\vskip 4mm

The present paper deals with the problem of 
characterizing the continuous functions on 
the unit circle which extend holomorphically 
into the unit disc, in terms of 
the argument principle. 
Generalizing some results of H.Alexander and 
J.Wermer [AW], E.L.Stout [S] 
obtained a characterization of continuous functions 
on boundaries of certain domains $D$ 
in $\C^n,\ n\geq 1$ which extend holomorphically 
through $D$, 
in terms of a generalized argument principle. In
the special case 
of $\D$, the open unit disc 
in $\C$, a version of his result is 
\vskip 3mm
\noindent\bf THEOREM 1.0\ \rm [S]\it\ A continuous
function $f$ on $b\D$ extends 
holomorphically through $\D$ if and only if
$$
\left.\eqalign
{&\hbox{for each polynomial\ }Q\hbox{\ of
two complex variables such}\cr &\hbox{that\ }
Q(z,f(z))\not= 0\ (z\in b\D)
\hbox{\ the change of argument of} \cr &\hbox{the 
function\ } z\mapsto Q(z,f(z))\hbox{\ around\ } b\D
\hbox{\ is nonnegative.}\cr}\ \ \ \ \right\}
\eqno (1.1)
$$
\rm J.Wermer [W] showed that it suffices to assume 
(1.1) only for polynomials of the form 
$Q(z,w)=w+P(z)$ provided that $f$ is smooth and asked 
whether the same holds for 
continuous functions. In the present paper we prove that 
this indeed is the case:
\vskip 3mm
\noindent\bf THEOREM 1.1\ \ \it A continuous function $f$ 
on $b\D$ extends 
holomorphically through $\D $ if and only if
$$
\left.\eqalign
{&\hbox{for each polynomial\ }P\hbox{\ such that\ }f+P\not=0\hbox{\ on\ }
b\D\ \hbox{\ the\ }\cr
&\hbox{change of argument of\ }f+
P\hbox{\ around\ }b\D \hbox{\ is nonnegative.\ }\cr}\ \ \ 
\right\}
\eqno(1.2)
$$
\rm Note that the only if part is an obvious 
consequence of the argument principle. In fact, if 
$f$ admits a holomorphic extension $\tilde f$ 
then the change of argument of $f+P$ around $b\D$ 
 equals $2\pi $ times the number of zeros of 
 $\tilde f + P$ in $\D$.
\vskip 6mm
\bf 2.\ The Morera condition\rm 
\vskip 4mm
\noindent \bf LEMMA 2.1\ \ \it Let $f$ be a continuous 
function on $b\D$ which satisfies 
\rm (1.2). \it Then
$$
\int_{b\D} f(z) dz = 0.
$$
\bf Proof.\ \rm Suppose that $\int_{b\D}f(z)dz \not= 0$. 
With no loss of generality assume that 
$$
{1\over{2\pi i}}\int_{b\D} f(z)dz = 1.
$$
Then $z\mapsto zf(z)-1$ is a continuous function on 
$b\D $ with zero average. Since 
the polynomials in $z$ and $\overline z$ are dense 
in $C(b\D)$ it follows that 
there are polynomials $P$ and $Q$ and a continuous 
function $g$ on $b\D$ such that 
$$
|g(z)|\leq 1/2\ \ (z\in b\D)
\eqno (2.1)
$$
and such that
$$
zf(z)-1= zP(z)+\overline{zQ(z)} + g(z)\ \ (z\in b\D).
$$
It follows that 
$$
z[f(z)-P(z)-Q(z)]\in 1+g(z)+i\R\ \ (z\in b\D)
$$
which, by (2.1), implies that
$$
{1\over 2} \leq \Im \bigl( z\bigl[ f(z)-P(z)-Q(z)\bigr]\bigr) 
\leq {3\over 2}\ \ (z\in b\D),
$$
so the change of argument of $z\mapsto z\bigl[ f(z)-P(z)-Q(z)\bigr]$ 
around $b\D$ equals zero. Thus, $f-P-Q\not=0$ on $b\D$ and the change 
of argument of $f-P-Q$ around $b\D$ equals $-2\pi $, 
contradicting the assumption that $f$ 
satisfies (1.2). This completes the proof.
\vskip 2mm
Before proceeding to the proof of Theorem 1.1 observe that Lemma 2.1 provides 
a simple alternative 
proof of Theorem 1.0; in fact, it provides a proof of the following corollary 
which sharpens Theorem 1.0. 
\vskip 3mm
\noindent \bf COROLLARY 2.1\ \it A continuous function $f$ on $b\D$ extends 
holomorphically 
through $\D$ if and only if for each nonnegative integer $n$ and each 
polynomial $P$ such that 
$z^nf(z)+P(z)\not= 0\ (z\in b\D )$ the change of argument of the function 
$z\mapsto z^nf(z)+P(z)$ 
around $b\D $ is nonnegative.
\vskip 2mm \noindent \bf Proof. \rm If $f$ satisfies the condition 
then Lemma 2.1 implies that
$$
\int_{b\D} z^n f(z) dz = 0
$$ for each nonnegative integer $n$. It 
is well known that this implies 
that $f$ extends holomorphically through $\D$. This completes 
the proof.
\vskip 6mm
\bf 3.\ Proof of Theorem 1.1 \rm
\vskip 4mm
\noindent\bf LEMMA 3.1\ \it Suppose that $f$ is a continuous 
function on $b\D $ that satisfies \rm (1.2)\it . Then for each $w\in\C\setminus\overline\D$ 
the function $z\mapsto f(z)/(z-w)$ satisfies \rm (1.2)\it . \rm
\vskip 2mm
\noindent\bf Proof. \rm Suppose that $w\in\C\setminus\overline\D$ and assume that $P$ is a
polynomial such that 
$$
{{f(z)}\over{z-w}}+P(z)\not = 0\ \ \ (z\in b\D ).
$$
Then
$$
f(z) + (z-w)P(z)\not = 0\ \ \ (z\in b\D)
$$
and since $f$ satisfies (1.2) it follows that the change of argument of the function 
$z\mapsto f(z)+(z-w)P(z)$ around $b\D $ is nonnegative. Since
$$
{{f(z)}\over{z-w}}+P(z) = {1\over{z-w}}\bigl[f(z)+(z-w)P(z)\bigr]\ \ (z\in b\D )
$$
and since the change of argument of $z\mapsto 1/(z-w)$ around $b\D$ is zero it follows that 
the change of argument of 
$$
z\mapsto {{f(z)}\over{z-w}}+P(z)
\eqno (3.1)
$$
around $b\D$ is nonnegative. This completes the proof.
\vskip 3mm
\noindent \bf Proof of Theorem 1.1. \rm Suppose that $f$ satisfies (1.2). By Lemma 3.1, for
each $w\in\C\setminus\overline\D$ the function (3.1) satisfies (1.2). By Lemma 
2.1 it follows that
$$
\int _{b\D } {{f(z)}\over{z-w}}dz = 0 \ \ (w\in\C\setminus\overline\D).
\eqno(3.2)
$$
It is well known that (3.2) implies that $f$ extends holomorphically 
through $\D$. This completes the proof. 
 \vskip 6mm
 \bf 4.\ An example \rm
\vskip 3mm
Theorem 1.1 gives a simple characterization of continuous 
functions on $b\D$ that extend 
holomorphically to $\D$ in terms of the argument principle. 
One can ask whether one can further simplify this 
characterization. J.Wermer [W] showed that in Theorem 1.1 it 
is not enough to assume  
(1.2) for 
polynomials of the form  
$P(z)=az+b, \ a,b,\in \C$. In this section we sharpen this by showing that for 
Theorem 1.1 to hold 
(1.2) must hold for polynomials of arbitrarily large degree. 
\vskip 3mm
\noindent\bf PROPOSITION 4.1\ \it For every $n_0\in\N$ there is a continuous 
function $f$ on $b\D$ such that 
whenever $P$ is a polynomial of degree not exceeding $n_0$ such that 
$f+P\not= 0$ on $b\D $, then the change of 
argument of $f+P$ around $b\D$ is nonnegative, yet $f$ does not extend 
holomorphically through 
$\D $.
\vskip 2mm
\noindent 
\bf Proof.\ \rm Let $n_0\in\N$ and let $n\in\N,\ n\geq n_0 + 1$. Set
$$
f(z) = z^n+{1\over{2z}}\ \ (z\in b\D).
$$
We show that for every polynomial $P$ of degree not exceeding $n_0$ such 
that $f+P\not=0$ on $b\D$, the change of argument of 
$f+P$ around $b\D$ is nonnegative. Assume, contrarily to what we want to 
prove, that there is a polynomial $P$, $\hbox{deg}(P)\leq n_0$, such 
that $f+P$ has no zero on $b\D$, and such that the change 
of argument of $f+P$ around 
$b\D$ is negative. 

Since $f+P$ is rational with a single pole in $\D$ the 
argument principle implies that 
the change of argument of $f+P$ around $b\D$ equals 
$2\pi (\nu -1)$ where $\nu $ is 
the number of zeros of $f+P$ on $\D $. By our assumption,
\ $2\pi (\nu-1) <0$, so 
$f+P$ has no zero on $\D$. The zeros of $f+P$ are the zeros of
$$
z^{n+1}+zP(z)+{1\over 2}.
\eqno (4.1)
$$
Since $n\geq n_0+1$ and since the degree of $P$ does 
not exceed $n_0$ it follows 
that the leading term in (4.1) is $z^{n+1}$. Since the 
constant term in (4.1) is 
$1/2$ it follows that the product of zeros of (4.1) 
equals 
$1/2$ which implies that at least 
one of the zeros of (4.1) is contained in $\D$, so 
at least one of the zeros of $f+P$ is 
contained in $\D$, a contradiction. This completes 
the proof.
\vskip 6mm
\bf 5.\ Holomorphic extendibility to finite Riemann surfaces \rm
\vskip 3mm
Theorem 1.0 has yet another, less elementary but 
even shorter proof using  
Wermer's maximality theorem: Suppose that $f$ is a 
continuous function on $b\D$ 
which satisfies (1.1) and which does not extend 
holomorphically through $\D$. 
By Wermer's maximality theorem [H] the polynomials in $z$ 
and $f$ are dense in $C(b\D )$. 
In particular, there is a polynomial $Q$ such that
$$
|Q(z,f(z))-\overline z|<{1\over 2} \ \ (z\in b\D).
$$
Obviously $Q(z,f(z))\not= 0 \ (z\in b\D)$ and the change 
of argument of 
$z\mapsto Q(z,f(z))$ around $b\D$ equals the change of 
argument of $z\mapsto\overline z$
 around $b\D$ which equals $-2\pi $. Thus, the change of 
 argument of $z\mapsto Q(z,f(z))$ 
 around $b\D$ is negative which contradicts the fact 
 that $f$ satisfies (1.1) and 
 so completes the proof of Theorem 1.1. This proof of 
 Stout's theorem was the first that the 
 author found. Only after a careful reading of Cohen's 
 proof of Wermer's maximality theorem 
 [C] the author found the proof of Lemma 2.1 which gives 
 a more elementary proof of
 Stout's theorem. J.Wermer has kindly pointed out to the 
 author that the preceding proof 
 generalizes to finite Riemann surfaces which gives 
 \vskip 3mm
 \noindent\bf THEOREM 5.1\ \it Let $X$ be a finite Riemann surface 
 with boundary $\Gamma $. Let 
 $A$ be the algebra of all continuous functions on $\Gamma $ which 
 extend holomorphically through $X$. A continuous function $f$ on $\Gamma $ 
 extends holomorphically through $X$ if and only if for 
 every polynomial $P$ with 
 coefficients in $A$ such that $P(f)\not= 0$ on $\Gamma $ 
 the change of argument of $P(f)$ 
 along $\Gamma $ is nonnegative.\rm
 \vskip 2mm
 \noindent\bf Proof.\ \rm Suppose that $f\in C(\Gamma )$ has the property 
 that the change of argument of 
 $P(f)$ along $\Gamma $ is nonnegative whenever $P$ is a polynomial with 
 coefficients in $A$ such that $P(f)
 \not= 0$ on $\Gamma $. Suppose that $f$ does not extend holomorphically 
 through $X$. By the maximality 
 of $A$ in $C(\Gamma )$ [R] the functions of the 
 form $P(f)$ where $P$ is 
 a polynomial with coefficients 
 in
 $A$, are dense in $C(\Gamma )$. Choose $g\in C(\Gamma )$
 such that $|g|=1$ on $\Gamma 
 $ and such that the
 change of argument of $g$ along $\Gamma $ equals $-2\pi $.
 There is a polynomial $P$ 
 with coefficients in $A$ such that  
 $|P(f)(z)-g(z)|<1/2\ (z\in\Gamma )$. Since $|g|=1$ on $\Gamma $
 the change of argument
 of $P(f)$ along $\Gamma $ equals 
 the change of argument of $g$ along $\Gamma $. 
 So the change of argument of $P(f)$ along $\Gamma $ is
 negative, contradicting the hypothesis. This proves that 
 $f$ extends holomorphically through $X$. The only if part 
 follows from the argument principle. This completes the proof.
 \vskip 4mm
  Mark Agranovsky observed that a substantially longer 
 argument in the proof of Theorem 1.1 in the original version of the paper can be replaced 
 by Lemma 3.1 which is due to him. The author is grateful to him for his kind permission to 
 include Lemma 3.1 into the final version of the paper. The author is also grateful to John 
 Wermer for pointing out that Theorem 5.1 follows from the maximality theorem. 
 
 The author is grateful to Larry Zalcman for the kind invitation 
 to include this paper into the proceedings of the Nahariya conference. 
 
 This work was supported 
 in part by the Ministry of Education, Science and Sport of Slovenia 
through research program Analysis and Geometry, Contract No.\ P1-0291. 
 \vskip 8mm
 \centerline{\bf REFERENCES}
\vskip 5mm
\noindent [AW]\ H.Alexander and J.Wermer: Linking 
numbers and boundaries of varieties. 

\noindent Ann.\ Math.\ 151 (2000) 125-150
\vskip 2mm
\noindent [C]\ P.J.Cohen: A note on constructive methods in Banach algebras.

\noindent Proc\ Amer.\ Math.\ Soc.\ 12 (1961) 159-163
\vskip 2mm
\noindent [H]\ K.\ Hoffman: Banach spaces of analytic functions.

\noindent Prentice-Hall, Englewood Cliffs, N.J., 1962
\vskip 2mm
\noindent [R]\ H.\ Royden: The boundary values of analytic and harmonic functions.

\noindent Math.\ Z.\ 78 (1962) 1-24
\vskip 2mm
\noindent [S]\ E.L.Stout: Boundary values and mapping degree.

\noindent Michig.\ Math.\ J.\ 47 (2000) 353-368
\vskip 2mm
\noindent [W] J.Wermer: The argument principle and boundaries of analytic varieties.

\noindent Oper. Theory Adv. Appl., 127, Birkhauser, Basel, 2001, 639-659
\vskip 15mm
\noindent Institute of Mathematics, Physics and Mechanics

\noindent University of Ljubljana,  Ljubljana, Slovenia

\noindent josip.globevnik@fmf.uni-lj.si

\bye